\documentclass{amsart}

\usepackage{amsmath,amssymb,amsbsy,amsfonts,amsthm,latexsym,amsopn,amstext, amsxtra,euscript,amscd} %, hyperref}

\usepackage{cite}

\usepackage{amsrefs}

\input xy \xyoption{all}

\newtheorem{thm}{Theorem}
\newtheorem{prop}{Proposition}
\newtheorem{lemma}{Lemma}
\newtheorem{remark}{Remark}

\newtheorem{exa}{Example}
\theoremstyle{definition}

\def\P{\mathbb P^1 (k)}
\def\Q{{\mathbb Q}}
\def\Z{{\mathbb Z}}
\def\C{{\mathbb C}}
\def\M{{\mathcal M}}

\def\L{\mathcal L}

\def\Aut{\mbox{Aut }}

\def\iso{\cong}

\def\D{\Delta}
\def\e{\varepsilon}

\def\g{\gamma}
\def\bG{\bar G}

\def\a{{\alpha }}

\def\<{\langle}
\def\>{\rangle}

\def\s{\mathfrak s}
\def\u{\mathfrak u}
\def\v{\mathfrak v}
\def\r{\mathfrak r}

\def\E{\mathfrak E}
\def\p{\mathfrak p}

\begin{document}

\title{Families of genus two curves with many elliptic subcovers.}

\author{Tony Shaska}

\address{Department of Mathematics and Statistics,   Oakland University, Rochester, MI, 48309.}

\email{shaska@oakland.edu}
\thanks{The author thanks ...}


\date{\today}


%    "Communicated by" -- provide editor's name; required.
\commby{xxx}

%****************
\begin{abstract} 
We determine all genus  2 curves, defined over $\C$,  which have simultaneously  degree 2 and 3  elliptic subcovers. The locus of such curves has three irreducible 1-dimensional genus zero  components in $\M_2$.  For each component we find a rational parametrization and  construct the equation of the corresponding genus 2 curve and its elliptic subcovers in terms of the parameterization.  Such families of genus 2 curves are determined for the first time.  Furthermore, we prove that there are only finitely many genus 2 curves (up to $\C$-isomorphism) defined over $\Q$,   which have degree 2 and 3 elliptic subcovers also defined over $\Q$. 

\end{abstract}

\maketitle

\section{Introduction}

If there is a degree $d$  covering from a genus 2 curve $C$ to an elliptic curve $E$ then the space $\L_d$ of such genus 2 curves is an algebraic 2-dimensional locus in $\M_2$ when $d \cong 1 \mod 2$.  Genus 2 curves with this property have been studied extensively in the XIX century.  In the last decade of the XX century there was renewed interests on the topic coming from interests from number theory, cryptography, mathematical physics, solitons, differential equations, etc.  These spaces for $d= 3, 5$ were explicitly computed by this author and his co-authors in \cite{deg3} and \cite{deg5}. 

Due to some of the interesting properties of such genus 2 curves they have found applications in cryptography, factoring of large numbers, etc.  There is always an interest in having genus 2 curves defined over $\Q$ with many elliptic subcovers.  

In \cite{Cosset} such genus two curves were used for factorization of large numbers. Although the arithmetic of $C$ is more complicated than on an elliptic curve, the author shows that this is balanced by the fact that each computation on $C$ essentially corresponds to a pair of computations carried out on the two elliptic curves $E_1$ and $E_2$.  

In this paper we give a family of genus 2 curves  which have 4 elliptic subcovers   such that two of them are of degree 2 and the other two of degree 3.  We determine such elliptic subcovers and the corresponding covers explicitly.  Let $\p=[C]$ denote the isomorphism class of a genus 2 curve, over $\C$, such that $\p \in \L_2 \cap \L_3$. We denote its degree 2 (resp. 3) elliptic subcovers by $E_1, E_2$  (resp. $E_3, E_4$).  Their $j$-invariants are denoted by $j_1, j_2, j_3, j_4$ respectively. 
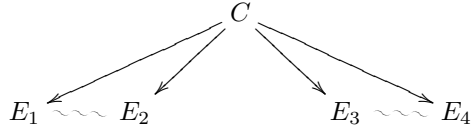
\[
%\begin{center}
\xymatrix{
 & & C  \ar@{->}[dl]  \ar@{->}[dll]   \ar@{->}[dr]  \ar@{->}[drr]    & & & \\
E_1  \ar@{~}[r]& E_2 & & E_3 \ar@{~}[r]& E_4 & \\
% & & &  E_3 \ar@{~}[r]& E_4 & \\
%\end{center}
}
\]

The locus $\L_2 \cap \L_3$ has three irreducible 1-dimensional components in $\M_2$.  Each component is a genus 0 curve and can therefore be parametrized. Using these parametrizations we are able construct 3 rational families of genus 2 curves and determine their 4 elliptic subcovers. Our main result can be summarized as below:

\begin{thm}[Main Theorem]\label{main_thm}
 a) The locus $\L_2 \cap \L_3$ of genus 2 curves which have four  elliptic subcovers, two of which are degree 2 and two of degree 3,  has three irreducible,   1-dimensional,  genus zero components in $\M_2$.\\
% which are determined explicitly in \eqref{locus1},  \eqref{locus2}, \eqref{locus3}. 

b) For all $t \in \C\setminus \{\D_t = 0\}$ there is a genus 2 curve $C_t$ where the $j$-invariants $j_1, j_2$ (resp. $j_3, j_4$)  of degree 2 (resp. degree 3) of elliptic subcovers are the roots of the quadratic $j^2+c_1 j+ c_0 =0$  (resp. given below): \\

i)  The equation of $C_t$ is
\begin{small}
\[ \begin{split}
y^2 = & 2\,{x}^{6}+ \left( 2\,t+12-2\,{t}^{3}+9\,{t}^{2} \right) {x}^{5}- \left( t+2 \right)  \left( {t}^{2}-2\,t-14 \right)  \left( {t}^{2}+1 \right) {x}^{4} - \\
& \left( {t}^{2}+1 \right)  \left( 2\,{t}^{4}-6\,{t}^{3}-31\,{t}^{2}-28\,t-37 \right) {x}^{3}- \left( {t}^{3}-10\,{t}^{2}-23\,t-28 \right)  \left( {t}^{2}+1 \right) ^{2}{x}^{2}\\
& + 6\, \left( t+2 \right)  \left( {t}^{2}+1 \right)^{3}x+2\, \left( {t}^{2}+1 \right)^{4},
\end{split}
\] 
\end{small}
where $\D_t= (t+2) (2t-11) (t^2+1) (t-1)$.   The $j$-invariants of elliptic subcovers are as follows,
\begin{small}
\[ \begin{split}
c_0 & =   \frac 1 {16} \,{\frac { \left( {t}^{2}-6\,t+4 \right)  \left( {t}^{2}+114\,t+ 124 \right)  \left( 2\,t-11 \right) ^{2}}{ \left( {t}^{2}+1 \right)^{4}}}\\   
c_1  &  = \,{\frac {-128} { \left( 2\,t-11 \right) ^{5}}}  \cdot  \left( 8\,{t}^{10}-320\,{t}^{9}+5580\,{t}^{8}-55460\,{t}^{7}+ 344593\,{t}^{6}-1379982\,{t}^{5}  \right.  \\
& \left.  +3562940\,{t}^{4}-4837160\,{t}^{3}+ 7580400\,{t}^{2}+180256\,t+1421824 \right)
\end{split} \]
\end{small}
\[ 
j_3  = 64\,{\frac { \left( {t}^{2}+114\,t+124 \right) ^{3}}{ \left( 2\,t-11
 \right) ^{5}}}, \qquad 
j_4  = 64\,{\frac { \left( {t}^{2}-6\,t+4 \right) ^{3}}{2\,t-11}}\\
\]

ii)  The equation of $C_t$ is  
\begin{small}
\[
\begin{split}
y^2    =  &  \left( 2916\,{x}^{3}t-2916\,{x}^{3}-486\,{x}^{2}{t}^{2}+1944\,{x}^{2}t-54\,x{t}^{4}+216\,x{t}^{3}-162\,x{t}^{2}+{t}^{7} \right.\\
& \left. -7\,{t}^{6}+15\,{t}^{5}-9\,{t}^{4} \right)  \left( 11664\,{x}^{3}+2916\,{x}^{2}-108\,x{t}^{3}+324\,x{t}^{2}+{t}^{6}-6\,{t}^{5}+9\,{t}^{4} \right)
\end{split}
\] 
\end{small}
where $\D_t = t (t-1) (t^2-2t-2) (t-3)$.   The $j$-invariants of elliptic subcovers are as follows, 
\begin{small}
\[ 
\begin{split}
c_0 & = - \frac 1 {16} \,{\frac { \left( -225+540\,t-396\,{t}^{2}+80\,{t}^{3} \right) 
 \left( t-3 \right) ^{6}}{{t}^{8} \left( t-1 \right) ^{2}}}
   \\   
c_1  &  =  {\frac {384} { \left( t-1
 \right)  \left( t-3 \right) ^{9}}} \cdot \left(-820125+4045950\,t+3221840\,{t}^{8}-5562270\,{t}^{2}- 33869934\,{t}^{5} \right. \\
 & \left. +24265899\,{t}^{4}-4786128\,{t}^{3}+25696944\,{t}^{6} -11732952\,{t}^{7}-491200\,{t}^{9}+32000\,{t}^{10}\right)  \\
\end{split}
\]
\end{small}
This family has isomorphic degree 3 elliptic subcover with $j$-invariants 
\[ 
j_3 \, = \, j_4  = -1728\,{\frac { \left( 2\,t-5 \right) ^{3}}{ \left( t-1 \right) ^{3}
 \left( t-3 \right) ^{3}}}
\]

iii) The equation of $C_t$ is  
\begin{small}  
\[\begin{split}
y^2  &  =      \left( 4\,{\frac { \left( 2\,{t}^{2}+7
\,t+2 \right) ^{6}{x}^{3}}{{t}^{2} \left( 2+3\,t+2\,{t}^{2} \right) ^{4}}}+{\frac { \left( 2\,{t}^{2}+7\,t+2 \right) ^{6}{x}^{2}}{{t}^{2}
 \left( 2+3\,t+2\,{t}^{2} \right) ^{4}}}+2\,{\frac { \left( 2\,{t}^{2}+7\,t+2 \right) ^{3}x}{t \left( 2+3\,t+2\,{t}^{2} \right) ^{2}}}+1
 \right)  \\
 & \left( {\frac { \left( 2\,{t}^{2}+7\,t+2 \right) ^{6}{x}^{3}}{{t}^{2} \left( 2+3\,t+2\,{t}^{2} \right) ^{4}}}-{\frac { \left( 2\,{t}^{2}+7\,t+2 \right)^{4} \left( 2\,{t}^{4}+{t}^{3}-5\,{t}^{2}-2\,t-1 \right) {x}^{2}}{{t}^{2} \left( 2+3\,t+2\,{t}^{2} \right)^{3}}}+{
\frac { \left( 2\,{t}^{2}+7\,t+2 \right)^{3}x}{t \left( 2+3\,t+2\,{t}^{2} \right)^{2}}}+1 \right) \\
\end{split}
\] 
\end{small}
where $\D_t=t (2t^2+7t+2) (t-2) (t+2) (2t-1) (t+1) (2+3t+2t^2)$.   The $j$-invariants of elliptic subcovers are as follows, 
\begin{small}
\[ 
\begin{split}
c_0 & =  1024\,{\frac { \left( t+2 \right) ^{6}{t}^{6} \left( {t}^{4}+56\,{t}^{
2}+16 \right)  \left( {t}^{4}-4\,{t}^{2}+16 \right) }{ \left( t-2
 \right) ^{2} \left( 2+3\,t+2\,{t}^{2} \right) ^{8}}}  \\   
c_1  &  = -\,  {\frac 4 {{t}^{4} \left( t-2 \right) ^{4} \left( t+2 \right) ^{4}}} 
\cdot \left( 4\,{t}^{16}-79\,{t}^{14}+1000\,{t}^{12}+3824\,{t}^{10}+ 207616\,{t}^{8} \right. \\
& \left. +61184\,{t}^{6}+256000\,{t}^{4}-323584\,{t}^{2}+262144 \right)   \\
\end{split}
\] \end{small}
\[ 
j_3:= 2\,{\frac { \left( {t}^{4}+120\,{t}^{3}+536\,{t}^{2}+480\,t+16 \right)^{3}}{t \left( t-2 \right) ^{8} \left( t+2 \right) ^{2}}}, \quad
j_4:= 256\,{\frac { \left( {t}^{4}-4\,{t}^{2}+1 \right) ^{3}}{{t}^{2} \left( t-2 \right)  \left( t+2 \right) }}
\] 

c) Every curve $[C]\in \L_2 \cap \L_3$ is isomorphic, over $\C$, to one of the curves in  i), ii), or iii) for some value of $t\in \C\setminus \D_t$. \\
\end{thm}

The rest of this paper will be proving this theorem. We will use the explicit equation of $\L_2$ and $\L_3$.  The idea of this paper is based on the following result where the explicit equation of $\L_3$ is computed. 

\begin{thm}[Shaska 2001]\label{deg3}
Let $K$ be a genus 2 field and $e_3(K)$ the number of $\Aut(K/k)$-classes of elliptic subfields of $K$ of degree 3.  Then;

i) $e_3(K) =0, 1, 2$, or $4$

ii) $e_3(K) \geq 1$ if and only if the classical invariants of $K$ satisfy  the irreducible  equation $F(J_2, J_4, J_6, J_{10})=0$ displayed in Appendix A in \cite{deg3}.
\end{thm}

The equation of the second part of the theorem is called the $\L_3$ equation throughout this paper. Its singularities were studied in \cite{beshaj}.  The equation of $\L_2$ has been known in different forms since the XIX century.  We will use the equation of $\L_2$ as in \cite{sh_2000} 

We use parametrizations of $\L_2$ by the $\s$-invariants which were introduced in \cite{sh_2000} and have been used by many authors since.  For computations in $\L_3$ we make use of the invariants among two cubics which were introduced in \cite{deg3} and seem to have been unknown to classical invariant theorists. 

After determining the locus $\L_2\cap \L_3$ in terms of absolute invariants $i_1, i_2, i_3$ of genus 2 curves we parametrize each component of this locus.  The constructing the genus 2 curves starting from the moduli point in these loci makes use of the Prop.~\ref{prop1} where the equation of the curve is given in terms of the $\s_1$, $\s_2$ invariants as in \eqref{generic_V4}.   Such equations, known to the author since 2003, are being published for the first time. For a method of how to determine a minimal equation of genus 2 curves over its field of definition check \cite{min_eq}.

%********************************************************************
\section{Preliminaries on genus 2 curves with split Jacobians}

Curves with split Jacobians have been studied extensively before by many authors.  For a survey on the topic and a complete list of references the reader can check \cite{sh_05}.  In this section we briefly set some of the notation and describe some results that we will need in the next section.  

\subsection{The space $\L_2$}

All our computations in this paper are based on the $\s$-invariants of genus 2 curves with extra involutions.  Hence we will define them here and describe some basic results.  For details the reader can check \cite{sh_05}. 

Let $C$ be a genus 2 curve and  $z_1$ is an elliptic involution of $C$. Denote by $\Gamma = PGL(2, \C)$. Let $z_2=z_1z_0$, where $z_0$ is the hyperelliptic
involution. Let $E_i$ be the fixed field of $z_i$ for $i=1,2$. 

We need to determine to what extent the normalization in the above proof determines the coordinate $X$. The condition $z_1(X)=-X$ determines the coordinate $X$ up to a coordinate change by some $\g\in \Gamma$ centralizing $z_1$. Such $\g$ satisfies $\g(X)=mX$ or $\g(X) = \frac m X$, $m\in k\setminus \{0\}$. 
Hence we can assume that the Weierstrass points are $\{ \pm \a_1, \pm \a_2, \pm \a_3\}$.  If we denote the symmetric polynomials of $\a_1^2, \a_2^2, \a_3^2$ by $a, b, c$ then $C$ has equation $Y^2= X^6 -a X^4 + bX^4 -c$.  The additional condition $abc=1$ forces $1=-\g(\a_1) \dots \g(a_6)$, hence $m^6=1$. Then $C$ is isomorphic to a curve with equation  
\begin{equation}
Y^2=X^6-aX^4+bX^2-1, 
\end{equation}
where  $27-18ab-a^2b^2+4a^3+4b^3\neq 0$.  

So $X$ is determined up to a coordinate change by the subgroup $H\iso D_6$ of $\Gamma$ generated by $\tau_1: X\to \e_6X$, $\tau_2: X\to
\frac 1 X$, where $\e_6$ is a primitive 6-th root of unity. Let $\e_3:=\e_6^2$. The coordinate change by
$\tau_1$ replaces $a$ by $\e_3b$ and $b$ by $\e_3^2 b$. The coordinate change by $\tau_2$ switches $a$ and $b$. Invariants of this action are:
\[ \s_1:=a b, \quad \s_2:=a^3+b^3 \]
The mapping
$$A: (\s_1,\s_2) \to (i_1, i_2, i_3)$$
gives a birational parametrization of $\L_2$. 

The ordered pair $\s_1, \s_2$ uniquely determines the isomorphism classes of curves in $\L_2$. 
%&&&&&&&&&&&&&&&&&&&&&&&&&&&&&&&&&&&&&&&&&
\begin{lemma} \label{lem1}
$k(\L_2) = k (\s_1, \s_2)$. 
\end{lemma}

The fibers of A of cardinality $>1$ correspond to those curves $C$ with $|\Aut(C)| > 4$.   The rational expressions  of $\s_1, \s_2$ can be found in \cite{sh_02, sh_05}

\subsection{Elliptic subcovers}
%&&&&&&&&&&&&&&&&&&&&&&&&&&&&&&&&&&&&&&&&&&&&&&&&&&&&&&&&&&&&

Let $j_1$ and $j_2$ denote the $j$-invariants of the elliptic curves $E_1$ and $E_2$. The
invariants $j_1$ and $j_2$ and are roots of the quadratic
\begin{equation} 
\begin{split}
j^2+256\frac {(2\s_1^3-54\s_1^2+9\s_1\s_2-\s_2^2+27\s_2)} {(\s_1^2+18\s_1-4\s_2-27)} j + 65536 \frac {(\s_1^2+9\s_1-3\s_2)} {(\s_1^2+18\s_1-4\s_2-27)^2} &
=0, \label{j_eq}
\end{split}
\end{equation}
see \cite{sh_2000} for details.

\begin{thm}
Let $\p = (\bar \s_1, \bar \s_2)  \in \L_2$ there exists a genus 2 curve $C_{\bar \s_1, \bar \s_2}$ with equation 
\begin{equation}\label{eq1}
Y^2 =a_0  X^6+a_1  X^5+a_2  X^4+ a_3 \, X^3+t \, a_2 X^2+t^2 a_1\, X +t^3  a_0, 
\end{equation}
where the coefficients are given by 
\begin{equation}\label{generic_V4}
\begin{split}
&    t =\bar \s_2^2-4 \bar \s_1^3 \\
&    a_0 =\bar \s_2^2+\bar \s_1^2 \bar \s_2-2 \bar \s_1^3  \\
&    a_1 =2 (\bar \s_1^2+3 \bar \s_2) \cdot (\bar \s_2^2-4 \bar \s_1^3)  \\
&    a_2 =(15 \bar \s_2^2-\bar \s_1^2 \bar \s_2-30 \bar \s_1^3) (\bar \s_2^2-4 \bar \s_1^3)\\
&    a_3 =4 (5 \bar \s_2-\bar \s_1^2)\cdot (\bar \s_2^2-4\bar \s_1^3)^2.\\
\end{split}   
\end{equation}
%
% and the expressions of $\s_1$ and $\s_2$ are given in terms of $i_1, i_2, i_3$ as in \cite{sh_02}. \\
\end{thm}

\proof  The proof can be computational.  Computing the absolute invariants $i_1, i_2, i_3$ we have 
\[
\begin{split}
i_1  & = \frac 9 4 \,{\frac {{\bar \s_1}^{2}-126\,\bar \s_1+405+12\,\bar \s_2}{ \left( 15+\bar \s_1 \right) ^{2}}}\\
i_2 & = {\frac {27}{8}}\,{\frac {729\,{\bar \s_1}^{2}+{\bar \s_1}^{3}+4131\,\bar \s_1-3645-1404\,\bar \s_2-36
\,\bar \s_1\bar \s_2}{ \left( 15+\bar \s_1 \right) ^{3}}} \\
i_3 & = {\frac {243}{8192}}\,{\frac { \left( -27-4\,\bar \s_2+18\,\bar \s_1+{\bar \s_1}^{2} \right) ^{
2}}{ \left( 15+\bar \s_1 \right) ^{5}}} \\
\end{split}
\]
Using the expressions of $\s_1, \s_2$ in \cite{sh_2000} in terms of $i_1, i_2, i_3$ we get \[ (\s_1, \s_2) = (\bar \s_1, \bar \s_2 ). \]  This completes the proof. 

\endproof 

That every genus 2 curve with automorphism group of order $> 2$ is defined over its field of moduli was proved before by Cardona/Quer and independently by this author in \cite{sh_02}.  
%There was a mistake in \cite{cardona} since the case when the Clebsch invariant $C=0$ was not treated.  
This expression of the curve in terms of the $\s$-invariants is the first one and it is beneficial because it uses the rational parametrization of the surface $\L_2$. 

We illustrate some of the ideas above with the following example. 

\begin{exa}
Let $C$ be a genus 2 curve with equation 
\[ \begin{split}
y^2  = \, & 3\,{x}^{6} +\left( 10\,\sqrt {3} - 8\right) {x}^{5} + \left( 63 -16\,\sqrt {3}\right){x}^{4}  + \left( 60\,\sqrt {3} -72\right) {x}^{3} \\
& +\left( 125 -40\,\sqrt {3}\right) {x}^{2} + \left( 46\,\sqrt {3} -36\right) x +29
\end{split} \]

Below we display some of the information about this curve from the Maple package "genus 2" written by the author.

\begin{verbatim}
Info(C,x);
   
    "The moduli point for this curve is p=(i_1, i_2, i_3) "
                                 [31  -125    361  ]
            (i[1], i[2], i[3]) = [--, ----, -------]
                                 [12   8    3981312]
 "The Automorphism group is isomorphic to the group with GapId"
                             [4, 2]
                      "Sh-invariants are "
                     (s[1], s[2]) = [3, 28]
                   "The  field of moduli is:"
                             M = Q
             "The minimal field of definition is:"
                             F = Q
                "The  degree of obstruction is:"
                         "[F : M]" = 1
  "Rational model is over its minimal field of definition is:"
        2        6          5            4              3
       y  = 491 x  + 62868 x  + 3615924 x  + 119727712 x 

                        2                               
          + 2444364624 x  + 28729167168 x + 151677646016

   "This curve has extra involutions.  Its degree 2 elliptic  subcovers 
   
   have j-invariants"
                    23584461610752  6750000
                    --------------, -------
                        312481       3913  
        "This curve has NO degree 3 elliptic subcovers"
\end{verbatim}

The rational model refereed above is computed using the  Eq.~\eqref{generic_V4}. By an appropriate M\"obius transformation one can show that it is isomorphic over $\C$ with the curve
\[ y^2= 491\,{x}^{6}+2418\,{x}^{5}+5349\,{x}^{4}+6812\,{x}^{3}+5349\,{x}^{2}+ 2418\,x+491 \]
It can be checked that it has the same $i_1, i_2, i_3$ invariants as the previous curve. In \cite{min_eq} we can show that we can do better 
and a more "minimal" equation.

\end{exa}

Hence, the equation provided in \eqref{generic_V4} is not necessary a "minimal" equation of the curve.  For a "minimal" equation of genus 2 curves see \cite{min_eq}. 

%****************************************************************
\subsection{The space $\L_3$}

In \cite{deg3} it was shown that every curve $C$ in $\L_3$ can be written as
\begin{equation}\label{eq_F1_F2}
Y^2=(X^3+a X^2+b X+1 )\, (4X^3+b^2X^2+2bX+1 )
\end{equation}
and the following
\begin{equation}\label{u_v}
\u \, = \, ab,  \quad \v   \, = \, b^3\\
\end{equation}
are invariants of $C$ under any change of coordinates.  

The mapping $k^2\setminus \{\D =0\}  \to \L_3$ such that 
\[ (\u, \v) \to (i_1, i_2, i_3) \]
has degree 2.  Instead the invariants of two cubics as defined in \cite{deg3}
\[
\begin{split}
\r_1 &= 27\frac {{\v}({\v}-9-2{\u})^3} {4{\v}^2-18{\u}{\v}+27{\v}-{\u}^2{\v}+4{\u}^3}\\
\r_2 & = -1296 \frac  {{\v}({\v}-9-2{\u})^4}   {({\v}-27)(4{\v}^2-18{\u}{\v}+27{\v}-{\u}^2{\v}+4{\u}^3)}
\end{split}
\]
uniquely determine the isomorphism class of curves in $\L_3$. However, for the curves in $\L_3$ the field of moduli is not necessary a field of definition.  One can show that the degree of the obstruction is   $\leq 2$ as proved below.

\begin{thm} $k(\r_1, \r_2) = k (\L_3)$.  Moreover, for every $\p = (\r_1, \r_2) \in \L_3$ there is a genus two curve $C$ with equation 
\begin{equation}\label{deg3} 
Y^2=(\v^2 X^3+\u \v X^2 +\v X+1)\, (4\v^2 X^3 +\v^2 X^2+2\v X+1),  
\end{equation}
\end{thm}

\proof
We compute absolute invariants $i_1, i_2, i_3$ in terms of $\u, \v$.  Substituting them in the equation of $\L_3$ we check that they satisfy this equation. 

\qed

We further discuss $\L_3$. We let
\[R:=(27 \v + 4 \v^2 - \u^2 \v + 4\u^3 -18 \u \v )\neq 0.\]
For $4\u - \v -9\neq 0$ the degree 3 coverings are given by  $\phi_1(X,Y)\to (U_1, V_1)$ and $\phi_2(X,Y)\to (U_2, V_2)$ where
\begin{small}
\begin{equation}\label{covers}
\begin{split}
& U_1= \frac {\v X^2} {\v^2 X^3 + \u \v X^2 + \v X +1}, \quad
U_2= \frac {(\v X +3)^2\, \,  (\v  (4\u-\v-9) X +3\u -\v)}
{\v \,(4\u -\v -9) (4 \v^2 X^3 + \v^2 X^2 + 2 \v X + 1) }, \\
&  \quad \quad
 V_1 =  Y \, \frac {\v^2 X^3 - \v X -2} {\v^2 X^3 + \u \v X^2 + \v X +1}, \\
& \quad \quad V_2 =   (27- \v)^{ \frac 3 2} \, Y
\frac {\v^2 (\v-4\u+8)X^3 +\v (\v-4\u)X^2 -\v X+1  }
{(4 \v^2 X^3 + \v^2 X^2 + 2 \v X + 1)^2 }\\
\end{split}
\end{equation}
\end{small}
and the elliptic curves have equations:
\begin{small}
\begin{equation}
\begin{split}
\E: & \quad V_1^2= R \, U_1^3 - (12 \u^2 - 2 \u \v - 18 \v )U_1^2 +
(12\u -\v) U_1 -4 \\
\E': & \quad V_2^2=c_3 U_2^3 +c_2 U_2^2 + c_1 U_2 +c_0
\end{split}
\end{equation}
\end{small}
where
\begin{small}
\begin{equation}
\begin{split}
c_0 & = - (9\u -2 \v -27 )^3\\
c_1 & = (4\u -\v -9) \, (729 \u^2 + 54 \u^2 \v -972\u \v - 18\u \v^2 +189\v^2 + 729 \v +\v^3) \\
c_2 & =- \v \,  (4\u -\v -9)^2 \,  (54\u +\u \v -27\v)     \\
c_3 & =\v^2 \, (4\u -\v -9)^3 \\
\end{split}
\end{equation}
\end{small}
The above facts can be deduced from Lemma 1 of \cite{deg3}. The case
$4\u - \v -9= 0$ is treated separately in \cite{deg3}.

There  is an   automorphism $\, \, \beta \in Gal_{k(\u,\v)/k(i_1, i_2, i_3)}$ given by
\begin{small}
\begin{equation}\label{eq_nu}
\begin{split}
\beta(\u) & =\frac {(\v-3\u)(324\u^2+15\u^2\v-378\u\v-4\u\v^2+243\v+72\v^2)}
{(\v-27)(4\u^3+27\v-18\u\v-\u^2\v+4\v^2)}\\
\beta(\v) & =- \frac {4(\v-3\u)^3}{4\u^3+27\v-18\u\v-\u^2\v+4\v^2} \\
\end{split}
\end{equation}
\end{small}
which permutes the $j$-invariants of $\E$ and $\E'$.
%The map $$\th: (\u,\v) \to (i_1, i_2, i_3)$$ defined when $J_2\neq 0$ and $\Delta \neq 0$   has degree 2.
These $j$ invariants are given explicitly in terms of $\u $ and $\v$ as below:
\begin{equation}\label{j_deg_3}
\begin{split}
j_3 & = 16\,{\frac {\v \left( 216\,{\u}^{2}+\v{\u}^{2}-126\,\u\v+405\,\v+12\,{\v}^{2}-
972\,\u \right) ^{3}}{ \left( 4\,{\u}^{3}-\v{\u}^{2}-18\,\u\v+4\,{\v}^{2}+27
\,\v \right) ^{2} \left( \v-27 \right) ^{3}}} \\
j_4 & = -256\,{\frac { \left( {\u}^{2}-3\,\v \right) ^{3}}{\v \left( 4\,{\u}^{3}-\v
{\u}^{2}-18\,\u\v+4\,{\v}^{2}+27\,\v \right) }} \\
\end{split}
\end{equation}
 
\begin{remark}
There are exactly two genus 2 curves (up to isomorphism) with $e_3(K)=4$, see 4.2 in \cite{deg3}. The case $e_3(K)=1$
(resp., 2) occurs for a 1-dimensional (resp., 2-dimensional) family of genus 2 curves.
\end{remark}

The theorem below shows that if we want to search for a family of curves with many elliptic subcovers we have to look at the curves with automorphism group $V_4$.

\begin{thm}[Shaska 2003]
Let $C$ be a genus two curve which has a degree 3 elliptic subcover. Then the automorphism group of $C$ is one of the
following: $\Z_2, V_4$, $ D_8$, or $D_{12}$. Moreover, there are exactly six curves $C\in \L_3$ with automorphism group $D_8$ and
six curves $C\in \L_3$ with automorphism group $D_{12}$.\\
\end{thm}

The list of all curves in $\L_3$ with automorphism group $> 4$ is given in \cite{sh_02}, where their rational points are determined also.

%For any prime $p \in \Z$ and a genus 2 curve $C$ we say that $C$ has a \textbf{good reduction} at $p$ if localized at $p$ it is still a genus 2 curve.  

\subsection{Constructing curves from their moduli points}

In our computations we will find the intersection $\L_2 \cap \L_3 $ as a sublocus in $\M_2$.  Hence, we need a constructive way to determine the equation of the curve once we know the moduli point.  We summarize all the results in the following. 

\begin{prop}\label{prop1} The following are true:
\begin{description}
\item [i)]  Let $j \in \Q$. Then there exists an elliptic curve $\E$ defined over $\Q$ such that $j(\E) = j$. 
Moreover the equation of $\E$ is given by \\

a) If $j \neq 0, 1728$ then  $ \Aut(\E)=\Z_2 $ and   \[ y^2  = x^3 - \frac j {48} (j-1728)^3 x - \frac j {864} (j-1728)^5\] 

b)  If $j=0$ then $\Aut (\E) = \Z_2 \times \Z_3$   and   $y^2= x^3 - \frac 1 4$. 

c) If  $j=1728$ then $\Aut (\E)=V_4$ and  $  y^2= x^3 + x$. \\

\item [ii)]  The space $\L_2$ is parametrized by the $\s$-invariants $(\s_1, \s_2)$, $k (\L_2) = k (\s_1, \s_2)$.  Let $\p \in \M_2$ such that $\p =  \left( i_1, i_2, i_3 \right) \in \Q^3$ and $\Aut (\p) \iso V_4$.  Then there exists a genus 2 curve $C$ defined over $\Q$ such that $\p = [C]$.  Moreover, its equation is 
\begin{equation}\label{eq1}
Y^2 =a_0  X^6+a_1  X^5+a_2  X^4+ a_3 \, X^3+t \, a_2 X^2+t^2 a_1\, X +t^3  a_0 
\end{equation}
where the coefficients are given by 
\begin{equation}\label{generic_V4}
\begin{split}
&    t =\s_2^2-4 \s_1^3 \\
&    a_0 =\s_2^2+\s_1^2 \s_2-2 \s_1^3  \\
&    a_1 =2 (\s_1^2+3 \s_2) \cdot (\s_2^2-4 \s_1^3)  \\
&    a_2 =(15 \s_2^2-\s_1^2 \s_2-30 \s_1^3) (\s_2^2-4 \s_1^3)\\
&    a_3 =4 (5 \s_2-\s_1^2)\cdot (\s_2^2-4\s_1^3)^2\\
\end{split}   
\end{equation}
and the expressions of $\s_1$ and $\s_2$ are given in terms of $i_1, i_2, i_3$ as in \cite{sh_02}. \\

\item [iii)] The space $\L_3$ is parametrized by the $\r_1, \r_2$-invariants as in \cite{deg3}, hence
$k(\L_3) = k(\r_1, \r_1)$.   Let $\p \in \L_3$.  Then there exists a genus 2 curve $C_\Q$ such that $p = [C]$ with equation
\begin{equation}\label{deg3} 
Y^2=(\v^2 X^3+\u \v X^2 +\v X+1)\, (4\v^2 X^3 +\v^2 X^2+2\v X+1). 
\end{equation}
\end{description}

\end{prop}

\proof The first part is elementary.  The proof of the part ii) can be found on \cite{sh_2000} or \cite{b-sh}. The equation of the curve $C$ given in Eq.\eqref{eq1} can be verified by computing the absolute invariants  of the curve and checking that they verify the equation of $\L_2$.   Since these computations are straight forward we do not display them here. For a rational moduli point $\p\in \M_2 (\Q)$,  the curve $C$ is defined over $\Q$ since $\s_1$ and $\s_2$ are given as rational functions in terms of $i_1, i_2, i_3$. 
A constructive proof of ii) and a discussion of a minimal polynomial of $C$ is intended in \cite{min_eq}. 

The proof of the third part iii) can be found in \cite{deg3}. 

\endproof

\begin{remark}
Invariants $(\s_1, \s_2)$ were also called $u, v$ in \cite{sh_02}.  Later they were generalized in \cite{g_sh} and used by many authors some of whom  have called them Shaska invariants or $\s$-invariants. We will call them $\s$-invariants not to confuse them with $\u, \v$ for degree 3 covers. 
\end{remark}

%It turns out that the equation in \eqref{eq1} is the minimal equation of $C$ over its minimal field of definition.  

%In the remaining of this section we explain the spaces $\L_2$ and $\L_3$ and how to obtain generic equation for a curve $\p = [C] $ in these spaces. 

%******************************************************************************************
\section{The locus $\L_2 \cap \L_3$.}

In order to construct genus 2 curves which have  degree 2 and degree 3 elliptic subcovers we need to determine the locus $\L_2\cap \L_3$.  This locus has three   components in $\M_2$, say 
\[ G_1 (i_1, i_2, i_3) \cdot G_2(i_1, i_2, i_3) \cdot G_3 (i_1, i_2, i_3) =0.\]  
We will show computational that each one of these components has genus zero. 
Parametrizing such components we are able to express the absolute invariants $i_1, i_2, i_3$ in terms of two variables $s$ and $t$ for all points $\p = (i_1, i_2, i_3) \in \L_2 \cap \L_3$.

Since for every point $\p \in \L_2$ the field of moduli is a field of definition then there is a curve $C$ with equations given as rational functions of $s$ and $t$ as in Prop.~\ref{prop1}, part ii).

\[ k^2\setminus \{ \D \neq 0\} \to \L_3 \cap \L_2 \to k^2\setminus \{ \D \neq 0\}\] 
\[ (\u, \v ) \to (i_1, i_2, i_3) \to (\s_1, \s_2) \]

The challenge here is to check the results of \cite{b-sh} and \cite{deg3} in order to see which ones are valid for curves and covers over $\Q$.  We know that for any rational point $\p \in \L_2 \cap \L_3$ the field of moduli is a field of definition.  In other words, there is a curve $C$ defined over $\Q$.  

%The map
%
%\[ (\u, \v ) \to \left( i_1, i_2, i_3   \right) \to \left( \s_1, \s_2 \right)\to (j_1, j_2) \]
%gives degree 2 elliptic subcovers.  

\subsection{$(u, v)$--space}   As it can be seen from above, it is a challenge computationally to lift from the point of moduli to the equation of the curve.  Instead we start with the curve given at Eq.~\eqref{deg3}.  Such curves are in $\L_2$ if and only if $\u$ and $\v$ satisfy the following 

\begin{small}
\begin{equation}\label{u_v_space}
\begin{split}
& (-18\, \u\, \v^2+\v^2\, \u^2+85\, \v^2-2160\, \v+468\, \u\, \v-28\, \v\, \u^2+4\, \u^3) \\
& (8\, \v^3+27\, \v^2-54\, \u\, \v^2-\v^2\, \u^2  +108\, \v\, \u^2+4\, \v\, \u^3-108\, \u^3)\\
& (3459375\, \v^3 -11390625\, \v^2 -333187\, \v^4+274410\, \u\, \v^3-1215000\, \u\, \v^2-324\, \v^6\\
& +2092500\, \v\, \u^3-1503225\, \v^2\, \u^2+374040\, \v\, \u^4-781106\, \v^2\, \u^3+443087\, \v^3\, \u^2\\
&  -69300\, \u\, \v^4+11168\, \v\, \u^5 -10864\, \v^2\, \u^4+24624\, \v^3\, \u^3-16535\, \u^2\, \v^4+2250\, \u\, \v^5\\
& +16929\, \v^5+128\, \u^6+81\, \v^5\, \u^2+54\, \v^4\, \u^3-16\, \v^4\, \u^4 +320\, \v^3\, \u^4 +32\, \v^3\, \u^5-1280\, \v^2\, \u^5)=0
\end{split}
\end{equation}
\end{small}

We can easily check that each  component has    genus 0.  Hence, we can parametrize each component. This parametrization will give   the equation of the curve $C$ and its degree 3 elliptic subcovers $\E$ and $\E^\prime$.  Computing such equations in each case will occupy the rest of this paper.

First we settle some notation. For a polynomial $F(x)= c_n x^n + \cdots + c_1 + c_0$ the \textbf{coefficient vector} we call the vector $(c_0, c_1, \dots c_n)^t$. 

The  parametrization methods used in some computational algebra packages as Maple, Mathematica etc might not produce the same results.  Indeed, in the third component their parametrizations were extremely long and we were not able to compute the equation of the genus 2 curve with such parametrizations.  
All our computations can be confirmed by substituting these parametrizations in each locus and verifying that the equation is satisfied. 
 
\subsection{First component} We start first with the component of the locus in Eq.~\eqref{u_v_space}, namely 
\[ -18 \u \v^2+\v^2 \u^2+85 \v^2-2160 \v+468 \u \v-28 \v \u^2+4 \u^3 =0. \]
This is a genus zero curve  which has a parametrization as follows
\[ \u = -(t+2)(t-4), \qquad \v = 2 \, \frac {(t+2)^3} {(t^2+1)} \]
Substituting these invariants in the expressions for  $\s_1, \s_2$ in \cite{sh_02, sh_05} we get
\begin{small}
\[ 
\begin{split}
\s_1 & = {\frac { \left( 4\,{t}^{2}-12\,t-5 \right)  \left( 8\,{t}^{2}+72\,t-13
 \right) }{ \left( 2\,t-11 \right) ^{2}}}   \\
\s_2 & = 
\frac 2       {\left( 2\,t-11 \right) ^{4}} \cdot
 \left(   128\,{t}^{8}-2560\,{t}^{7}+19776\,{t}^{6}-61248\,{t}^{5}+153600\,{t}^{4} \right. \\
 & \left. +185856\,{t}^{3}-192040\,{t}^{2}-33448\,t-8661  \right) \\
\end{split}
\]
\end{small}
The elliptic subcovers of degree 2 are determined by Eq.~\eqref{j_eq}. 
The $j$-invariants of degree 3 elliptic subcovers are obtained by replacing for $\u$ and $\v$ in Eq.~\eqref{j_deg_3}. 
\[ j_3=  64\,{\frac { \left( {t}^{2}+114\,t+124 \right) ^{3}}{ \left( 2\,t-11  \right) ^{5}}}, 
\quad j_4 = 64\,{\frac { \left( {t}^{2}-6\,t+4 \right) ^{3}}{2\,t-11}}\]
It is easy now to compute the equation of $C$, which has equation $Y^2=F(X)$ where $F(X)$ has coefficient vector
\[ 
 \left[ \begin {array}{c} 1\\ \noalign{\medskip}6\,{\frac { \left( t+2
 \right) ^{3}}{{t}^{2}+1}}\\ \noalign{\medskip}12\,{\frac { \left( t+2
 \right) ^{6}}{ \left( {t}^{2}+1 \right) ^{2}}}-2\,{\frac { \left( t+2
 \right) ^{4} \left( t-4 \right) }{{t}^{2}+1}}\\ \noalign{\medskip}20
\,{\frac { \left( t+2 \right) ^{6}}{ \left( {t}^{2}+1 \right) ^{2}}}+8
\,{\frac { \left( t+2 \right) ^{9}}{ \left( {t}^{2}+1 \right) ^{3}}}-8
\,{\frac { \left( t+2 \right) ^{7} \left( t-4 \right) }{ \left( {t}^{2
}+1 \right) ^{2}}}\\ \noalign{\medskip}48\,{\frac { \left( t+2
 \right) ^{9}}{ \left( {t}^{2}+1 \right) ^{3}}}-8\,{\frac { \left( t+2
 \right) ^{10} \left( t-4 \right) }{ \left( {t}^{2}+1 \right) ^{3}}}
\\ \noalign{\medskip}-32\,{\frac { \left( t+2 \right) ^{10} \left( t-4
 \right) }{ \left( {t}^{2}+1 \right) ^{3}}}+16\,{\frac { \left( t+2
 \right) ^{12}}{ \left( {t}^{2}+1 \right) ^{4}}}\\ \noalign{\medskip}
64\,{\frac { \left( t+2 \right) ^{12}}{ \left( {t}^{2}+1 \right) ^{4}}
}\end {array} \right] 
\]
The coefficients of $F(X)$ are obtained by simply replacing $s_1, s_2$ in Eq.~\eqref{generic_V4}.  

This completes the proof of the Main Theorem, part b), i).  

%************************************************************************************************
\subsection{Second component}  The second component 
\[ 8\, \v^3+27\, \v^2-54\, \u\, \v^2-\v^2\, \u^2  +108\, \v\, \u^2+4\, \v\, \u^3-108\, \u^3 =0, \] 
of Eq.~\eqref{u_v_space} is a genus zero curve.  We find a parametrization of this curve  as follows
\[ \u=9\,{\frac {t-4}{t \left( t-3 \right)  \left( t-1 \right) }}, \quad \v=-54\,{\frac {1}{{t}^{2} \left( t-3 \right) }}\]
Substituting these expressions into $\s_1, \s_2$ we get
\begin{small}
\[
\begin{split}
\s_1 & = -3 \, { \frac { -243+324\,t + 318\,{t}^{2}- 540\,{t}^{3}+ 125\,{t}^{4} }   { \left( t-3 \right)^4 }  } \\
\s_2 & =   \frac 6  {\left( t-3 \right)^8  }  \left( 59049-157464\,t+8748\,{t}^{2}+320760\,{t}^{3}-305802\,{t}^{4} \right.  \\
& \left.  +7128\,{t}^{5}+114540\,{t}^{6}-54200\,{t}^{7}+7625\,{t}^{8}  \right) \\
 \end{split}
 \]
\end{small} 
The elliptic subcovers of degree 2 are determined by Eq.~\eqref{j_eq}. 
The $j$-invariants of degree 3 elliptic subcovers are obtained by replacing for $\u$ and $\v$ in Eq.~\eqref{j_deg_3}.
Then $\E$ and $\E^\prime$ are isomorphic to each other and that is not very interesting to us, since we are looking for families with as many as possible elliptic subcovers. The j-invariant of such curves are 
\[ j_1=j_2= -1728\,{\frac { \left( 2\,t-5 \right) ^{3}}{ \left( t-1 \right) ^{3}  \left( t-3 \right) ^{3}}}\]

The genus 2 curve has equation 
\[
\begin{split}
 y^2 & =  ( 2916\,t{x}^{3}-2916\,{x}^{3}-486\,{t}^{2}{x}^{2}+1944\,t{x}^{2}-54\,{t}^{4}x+216\,{t}^{3}x-162\,{t}^{2}x   \\
 & +{t}^{7}-7\,{t}^{6}+15\,{t}^{5}-9\,{t}^{4}  )   \, ( 11664\,{x}^{3}+2916\,{x}^{2}-108\,{t}^{3}x+324\,{t}^{2}x+{t}^{6}\\
 &-6\,{t}^{5}+9\,{t}^{4}  ) 
 \end{split}
 \]

This completes the proof of the Main Theorem, part b), ii).  

%*********************************************************************
\subsection{Third component}  Next we consider  the third component of the locus in Eq.~\eqref{u_v_space}. 
We get a parametrization 
\[ 
\begin{split}
\u & =   -{\frac { \left( 2\,{t}^{2}+7\,t+2 \right)  \left( 2\,{t}^{4}+{t}^{3}-
5\,{t}^{2}-2\,t-1 \right) }{t \left( 2+3\,t+2\,{t}^{2} \right) }}\\
\v  & = {\frac { \left( 2\,{t}^{2}+7\,t+2 \right) ^{3}}{t \left( 2+3\,t+2\,{t}
^{2} \right) ^{2}}}\\
\end{split}
\]
which can be verified by substituting these expressions for $\u$ and $\v$ in the corresponding locus.  

Then  $\s_1, \s_2$ in terms of the parameter $t$ are as follows
\begin{small}
\[ 
\begin{split}
\s_1 & = {\frac { \left( {t}^{4}+3\,{t}^{3}+2\,{t}^{2}+6\,t+4 \right)  \left( 4 \,{t}^{4}+9\,{t}^{3}+26\,{t}^{2}+12\,t-8 \right) }{{t}^{3} \left( t-2
 \right)  \left( t+2 \right) ^{3}}}  \\
\s_2 & = \frac { 1 }   {4\,{t}^{5} \left( t-2 \right) ^{2} \left( t+2 \right) ^{5}} \cdot ( 16\,{t}^{14}+208\,{t}^{13}+896\,{t}^{12}+2940\,{t}^{11}+7785\,{t}^{10}\\
&+16926\,{t}^{9}+22832\,{t}^{8}+18272\,{t}^{7}-4640\,{t}^{6}-42816\,{t}^{5}-50688\,{t}^{4}-35328\,{t}^{3}\\
& -30464\,{t}^{2}-18944\,t-4096
 )  \\
\end{split}
\]
\end{small}

The $j$-invariants of the degree 3 elliptic subcovers are 

\[
\begin{split}
\E : \quad j_1 & = 2\,{\frac { \left( {t}^{4}+120\,{t}^{3}+536\,{t}^{2}+480\,t+16
 \right) ^{3}}{t \left( t-2 \right) ^{8} \left( t+2 \right) ^{2}}}
\\
\E^\prime : \quad j_2 & = 256\,{\frac { \left( {t}^{4}-4\,{t}^{2}+1 \right) ^{3}}{{t}^{2}
 \left( t-2 \right)  \left( t+2 \right) }}\\
\end{split}
\]

Now we can compute the equation of $C$ via Eq.~\eqref{generic_V4} and we get $Y^2=F(X)$ where $F(X)$ has coefficient vector 
\[
\left[ \begin {array}{c} 1\\ \noalign{\medskip}3\,{\frac { \left( 2\,
{t}^{2}+7\,t+2 \right) ^{3}}{t \left( 2+3\,t+2\,{t}^{2} \right) ^{2}}}
\\ \noalign{\medskip}3\,{\frac { \left( 2\,{t}^{2}+7\,t+2 \right) ^{6}
}{{t}^{2} \left( 2+3\,t+2\,{t}^{2} \right) ^{4}}}-{\frac { \left( 2\,{
t}^{2}+7\,t+2 \right) ^{4} \left( 2\,{t}^{4}+{t}^{3}-5\,{t}^{2}-2\,t-1
 \right) }{{t}^{2} \left( 2+3\,t+2\,{t}^{2} \right) ^{3}}}
\\ \noalign{\medskip}5\,{\frac { \left( 2\,{t}^{2}+7\,t+2 \right) ^{6}
}{{t}^{2} \left( 2+3\,t+2\,{t}^{2} \right) ^{4}}}+{\frac { \left( 2\,{
t}^{2}+7\,t+2 \right) ^{9}}{{t}^{3} \left( 2+3\,t+2\,{t}^{2} \right) ^
{6}}}-2\,{\frac { \left( 2\,{t}^{2}+7\,t+2 \right) ^{7} \left( 2\,{t}^
{4}+{t}^{3}-5\,{t}^{2}-2\,t-1 \right) }{{t}^{3} \left( 2+3\,t+2\,{t}^{
2} \right) ^{5}}}\\ \noalign{\medskip}6\,{\frac { \left( 2\,{t}^{2}+7
\,t+2 \right) ^{9}}{{t}^{3} \left( 2+3\,t+2\,{t}^{2} \right) ^{6}}}-{
\frac { \left( 2\,{t}^{2}+7\,t+2 \right) ^{10} \left( 2\,{t}^{4}+{t}^{
3}-5\,{t}^{2}-2\,t-1 \right) }{{t}^{4} \left( 2+3\,t+2\,{t}^{2}
 \right) ^{7}}}\\ \noalign{\medskip}-4\,{\frac { \left( 2\,{t}^{2}+7\,
t+2 \right) ^{10} \left( 2\,{t}^{4}+{t}^{3}-5\,{t}^{2}-2\,t-1 \right) 
}{{t}^{4} \left( 2+3\,t+2\,{t}^{2} \right) ^{7}}}+{\frac { \left( 2\,{
t}^{2}+7\,t+2 \right) ^{12}}{{t}^{4} \left( 2+3\,t+2\,{t}^{2} \right) 
^{8}}}\\ \noalign{\medskip}4\,{\frac { \left( 2\,{t}^{2}+7\,t+2
 \right) ^{12}}{{t}^{4} \left( 2+3\,t+2\,{t}^{2} \right) ^{8}}}
\end {array} \right] 
\]
  
The reader can easily verify that this equation of the genus 2 curve is the same with that claimed in the Main Theorem. This completes the proof of the Main Theorem, part b), i).  \\

\medskip  
  
\noindent \textbf{Proof of the Main theorem:}    Combining computations for each component   we have the results of Thm.~\ref{main_thm}.  

\qed

Given any genus 2 curve  $C$ over the complex numbers, we can check if it has degree 2 and degree 3 elliptic subcovers.  Computationally this is possible from our \textbf{genus2} Maple package. 
%For details see \cite{genus2}.  
Below we illustrate with an example.

\begin{exa}
Let $C$ be a genus 2 curve with equation
\[ y^2 = 8\,{x}^{6}+187\,{x}^{4}-1355\,{x}^{2}-1088\,{x}^{3}+3993\,x-3\,{x}^{5}
+2730
 \]
In the \textbf{genus2} package we enter only the polynomial $f(x)$. 

\begin{tiny}
\begin{verbatim}
        6        4         2         3               5       
f := 8 x  + 187 x  - 1355 x  - 1088 x  + 3993 x - 3 x  + 2730:
Info( f, x);
                "Initial equation of the curve"
  2      6        4         2         3               5       
 y  = 8 x  + 187 x  - 1355 x  - 1088 x  + 3993 x - 3 x  + 2730
"Igusa invariants are [J_2, J_4, J_6, J_10]"
       [5435864, 11848141844056, -15448925968029277668, -954225510546747438407778509664]
"Clebcsh invariants are [A, B, C, D]"
        [-679483  11104643929484  -1875304637846687089  3694376275475287252631213749681232]
        [-------, --------------, --------------------, ----------------------------------]
        [  15          5625              15625                      7119140625            ]

"The moduli point for this curve is p=(i_1, i_2, i_3) "
                      [4624326  -54010420569  -90812685325761]
 (i[1], i[2], i[3]) = [-------, ------------, ---------------]
                      [ 80089     45330374    929398866761216]
"The Automorphism group is isomorphic to the group with GapId"
                             [4, 2]
"Sh-invariants are "
                                [-688  -117101]
                 (s[1], s[2]) = [----, -------]
                                [ 27     972  ]
"The degree 2 j-invariants are roots of the quadratic"

"The  field of moduli is:"
                             M = Q
"The minimal field of definition is:"
                             F = Q
"The  degree of obstruction is:"
                         "[F : M]" = 1
"Rational model is over its minimal field of definition is:"

"A minimal rational model is over its minimal field of definition is:" 
  2                                6                                    5                                      4
y = -73190944927004215253422952448x + 111056747323558270032344844748800x +152821713034110852986785299725241600x  
                                          3                                             2
- 77928924518649965576021241379266720000 x + 12331972128546830355070565201074691490000 x 
+ 723169150606597090269353486513192776125000 x- 38459157573904353552637009599366812535828125
"with moduli point"
            [4624326  -54010420569  -90812685325761]
            [-------, ------------, ---------------]
            [ 80089     45330374    929398866761216]
"This curve has degree 3 elliptic subcovers.  "
\end{verbatim}
\end{tiny}
\end{exa}

\begin{thm}
There are only finitely many genus 2 curves (up to isomorphism) defined over $\Q$ such that they have degree 2 and degree 3 elliptic subcovers also defined over $\Q$.
\end{thm}

\proof  Given a genus two curve $C $ in the locus $\L_2 \cap \L_3$, it is defined over $\Q$ is and only if the moduli point $\p = (i_1, i_2, i_3)$ is a rational point.  

Let $\p =(i_1, i_2, i_3) \in \L_2\cap \L_3$.  Then, this $\p$ is in one of the components from Thm.~\ref{main_thm}. Since the corresponding $\s$-invariants are defined as rational functions in terms of $i_1, i_2, i_3$, then $\s_1, \s_2 \in \Q$.  Moreover, $C$ is defined over $\Q$ if and only if $\s_1, \s_2 \in \Q$. 

The invariants $j_1, j_2$ are rational numbers if 
\[ 
\begin{split}
g(\s_1, \s_2) &:=  -2916\,{s_{{1}}}^{2}s_{{2}}-864\,{s_{{1}}}^{3}s_{{2}}+486\,s_{{1}}{s_{{2}}}^{2}+189\,{s_{{2}}}^{2}{s_{{1}}}^{2}+2916\,{s_{{1}}}^{4}-216\,{s_{{1}}}^{5}\\
& -54\,{s_{{2}}}^{3}+729\,{s_{{2}}}^{2}+36\,{s_{{1}}}^{4}s_{{2}}-4\,{s_{{1}}}^{3}{s_{{2}}}^{2}-18\,s_{{1}}{s_{{2}}}^{3}+4\,{s_{{1}}}^{6}+{s_{{2}}}^{4}-4\,{s_{{1}}}^{2}\\
& -36\,s_{{1}}+12\,s_{{2}}
\end{split}
 \]
is a complete square in $\Q$, where $g(\s_1, \s_2)$ is the discriminant of the quadratic in \ref{main_thm}.  This expression is a complete square in $\Q$ if and only if the the curve 
\[ z^2 = g(\s_1 , \s_2) \]
has rational points.  
The curve $z^2 = g(\s_1 , \s_2)$ has genus 13, 13, 23 when $\p$ is in the locus $G_1, G_2, G_3$ respectively. From Falting's theorem, it has only finitely many rational points. This completes the proof. 

\endproof
%*********************************************

\bibliographystyle{amsplain}

\end{document}